\documentclass[10pt]{amsart}
\usepackage[all]{xy}
\usepackage{amssymb}

\newcommand{\Spe}{{\rm Spec}}

\newcommand{\Dgl}{{D_{\rm gl}}}

\newcommand{\TO}{{\mathcal{T}_\mathcal{O} }}
\newcommand{\T}{{\mathcal{T}}}
\newcommand{\A}{{ \rm Aut }}

\newcommand{\F}{{\mathbb{F}}}

\renewcommand{\O}{{\mathcal{O}}}
\newcommand{\Z}{{\mathbb{Z}}}
\newcommand{\lc}{\left\lceil}
\newcommand{\rc}{\right\rceil}
\newcommand{\lf}{\left\lfloor}
\newcommand{\rf}{\right\rfloor}

\renewcommand{\mod}{{\;\rm mod}}

\newtheorem{pro}{Proposition}[section]
\newtheorem{lemma}[pro]{Lemma}
\newtheorem{rem}[pro]{Remark}
\newtheorem{cor}[pro]{Corollary}

\newtheorem{orism}[pro]{Definition}

\begin{document}
\bibliographystyle{amsplain}

\title{Deformation of Curves with automorphisms and representations on Riemann-Roch spaces. }

\author{A. Kontogeorgis}

\begin{abstract}
We study the deformation theory of nonsigular projective curves 
defined over algebraic closed fields of positive characteristic.
We show that under some assumptions  
the local deformation problem for automorphisms of powerseries 
can be reduced to a deformation problem for matrix representations.
We study both equicharacteristic and mixed deformations 
in the case of two dimensional representations.
\end{abstract}

\address{Max-Planck-Institut f\"ur Mathematik
Vivatsgasse 7
53111 Bonn, Germany  
\newline
and 
\newline
$\;\;\;$ 
Department of Mathematics, University of the \AE gean, 83200 Karlovassi, Samos,
Greece\newline { \texttt{\upshape http://myria.math.aegean.gr/\~\;kontogar}}
}

\email{kontogeo@mpim-bonn.mpg.de,kontogar@aegean.gr}

\date{\today}

\maketitle
\section{Introduction}
Let $X$ be a nonsingular projective curve  of genus $g\geq 2$ defined over an algebraically closed field 
of characteristic $p>0$. The automorphism group $G:=\A(X)$ is known to be  a finite group.
The appearance of wild ramification in the cover $X \rightarrow X/\A(X)$ makes the 
theory of such covers more difficult than the corresponding theory in characteristic zero.
For a point $P \in X $ the decomposition group $G(P)=\{\sigma\in G: \sigma(P)=P\}$ is known to be cyclic in characteristic zero 
and a non-abelian solvable group admitting a ramification filtration \cite{SeL}.
In \cite{KontoMathZ} the author defined a faithful representation of the $p$-part of the decomposition 
group at a wild ramified point $P$:

\begin{equation} \label{fai-rep}
\rho: G_1(P) \rightarrow GL( L(mP)),
\end{equation}
 where $L(mP)=\{f \in k(X): \mathrm{div}(f) +mP \geq 0\} \cup \{ 0\}$. 
In this paper we would like to study the relation of two deformation theories, namely 
the deformation theory of representations of finite groups and the deformation theory of curves with automorphisms.

We will treat both  mixed characteristic and 
equicharacteristic deformations.
For the mixed characteristic case we consider $\Lambda$ to be a complete Noetherian local ring with residue field 
$k$. Usually $\Lambda$ is an algebraic extension of the ring of Witt 
vector $W(k)$.  For the equicharacteristic case we take $\Lambda=k$.

Let $\mathcal{C}$ denote the category of local Artin $\Lambda$-algebras,
which residue field $k$.
Consider    a subgroup $G$ of the group  $\A(X)$.
A deformation of the couple $(X,G)$ over the local Artin ring $A$ is a proper, smooth family
of curves 
\[
\mathcal{X} \rightarrow  \Spe(A)
\]
parametrized by the base scheme $\Spe(A)$, together with a group homomorphism $G\rightarrow \A _A(\mathcal{X})$, such that there is a 
$G$-equivariant isomorphism $\phi$ 
from 
the fibre over the closed 
point of $A$ to the original curve $X$:
\[
\phi: \mathcal{X}\otimes_{\Spe(A)} \Spe(k)\rightarrow X. 
\]
Two deformations $\mathcal{X}_1,\mathcal{X}_2$ are considered 
to be equivalent if there is a $G$-equivariant isomorphism $\psi$ that 
reduces to the identity in the special fibre and 
making the following diagram  commutative:
\[
\xymatrix{
\mathcal{X}_1 \ar[rr]^{\psi} \ar[dr] & & \mathcal{X}_2 \ar[dl] \\
& \Spe A &
}
\]
The  global deformation functor is defined:
\[
\Dgl: \mathcal{C} \rightarrow \rm{Sets}, 
A \mapsto
\left\{
\mbox{
\begin{tabular}{l}
Equivalence classes \\
of deformations of \\
couples $(X,G)$ over $A$
\end{tabular}
}
\right\}
\]
By the local-global theorems of J.Bertin and A. M\'ezard  \cite{Be-Me} and the formal patching theorems of 
D. Harbater,  K. Stevenson \cite{HarMSRI03}, \cite{HarStevJA99},  the study of the functor $\Dgl$ can be reduced to the study of the following deformation 
functors  attached to each wild ramification point  $P$ of the cover $X \rightarrow X/G$:
\begin{equation} \label{Bertin-Mezard-functor}
D_P:\mathcal{C} \rightarrow {\rm Sets}, 
A \mapsto 
\left\{
\mbox{
\begin{tabular}{l}
lifts $G(P)\rightarrow \A (A[[t]])$ of $\rho$ mod- \\ulo 
conjugation with an element \\ of $\ker(\A A[[t]]\rightarrow k[[t]] )$
\end{tabular}
}
\right\}
\end{equation} 
The theory  of automorphisms of formal powerseries rings is not as well understood as  is 
the theory of automorphisms of finite dimensional vector spaces, i.e the theory of general linear groups.

For a  $k$-algebra $A$ with maximal ideal $m_A$, 
consider the multiplicative group $L_n(A)<GL_n(A)$, of invertible  
lower  triangular matrices with entries in $A$, and invertible 
elements $\lambda$ in the diagonal, such that $\lambda-1 \in m_A$.
We consider the following  functor from the category $\mathcal{C}$ 
of local Artin $k$-algebras to the category of sets 
\begin{equation} \label{Fdeformation}
F: A \in Ob(\mathcal{C}) \mapsto \left\{
\begin{array}{l}
\mbox{liftings of } \rho: G(P) \rightarrow L_n(k) \\
\mbox{to } \rho_A: G(P) \rightarrow L_n(A) 
\mbox{ modulo} \\ \mbox{conjugation by an element }\\
\mbox{of }  
\ker(L_n(A)\rightarrow L_n(k))
\end{array}
\right\}
\end{equation}

It is known that among the curves  $X$ with automorphism group $G=\A(X)$ divisible by the 
characteristic, the curves so that $G_2(P)=\{1\}$ for all ramified points are the most simple. 
We will call these curves {\em weakly ramified}.
Many intractable problems for the theory of curves with general  automorphism group are solved 
for weakly ramified curves. For example  the computation of the $G$-module structure of spaces 
of holomorphic differentials  \cite{Koeck:04} or the computation of the deformation rings of curves with automorphisms \cite{CK}.
In our representation perspective it seems that the simplest curves are those with  two 
dimensional representations at all wild ramified points.
Notice that if  a  two dimensional representation 
 is attached at the wild point $P$,
 then the group $G_1(P)$ is elementary abelian 
and has conductor $m>1$ \cite[example 3.]{KontoMathZ}.

In section \ref{localMAT} we show how to attach a deformation 
of a matrix representation to every 
deformation of the couple $(X,G)$ over a complete local domain.
Section \ref{MatrixDeformations} is devoted to the deformations 
of matrix representations. We focus on the two dimensional 
case and we construct a hull for these deformations.
The deformation theory of such representations is
closely related to deformations of products of 
$\mathbb{G}_a$ group schemes  in the equal characteristic case 
or to $\mathcal{G}^{(\lambda)}$ group  schemes 
in mixed characteristic. In section \ref{SmallExt}
we try to analyze further the relation 
between the functors $F(\cdot)$ and $D(\cdot)$. 
A matrix representation allows us to 
express a deformation $\tilde{\rho}_\sigma$ given 
as a formal series  $\tilde{\rho}_\sigma(t) \in A[[t]]$
in the form of a root od rational function of $t$.
For the case of two dimensional representations,
where $V=G_1(P)$ is an 
elementary abelian group, we are able to 
compute the image of 
elements in $F(\cdot)$ in the tangent space $D(k[\epsilon]/\epsilon^2)=
H^1(V,\TO)$, see proposition \ref{5.1}.
 By combining these results to 
the computation of $H^1(V,\TO)$ given by the author in 
\cite[prop. 2.8]{KontoANT} we are able to compute the 
Krull dimension of the hull's attached to every wild 
ramified point. 

Finally, in section \ref{PRIES} we restrict ourselves to 
 to the equicharacteristic case and we 
 relate two dimensional matrix deformations to the deformation 
functor of R. Pries.

{\bf Acknowledgments} The author would like  to thank the participants of the conference
in Leiden on 
{\em Automorphisms of Curves} for enlightening 
conversations and especially R. Pries and M. Matignon for their corrections and remarks. This paper was completed during the author's visit at Max-Planck 
Institut f\"ur Mathematik in Bonn. The author would like to thank 
this insitution for its support and hospitality.

\section{Branch locus and liftings of matrix representations.}
\label{localMAT}
In this section we will show how the problem of deforming 
the representations attached at the wild ramified points give information 
on the problem of deformations of curves with automorphisms.

Select a wild ramified point $P_i$ on every orbit of  wild ramified points under the action of the group $G$.
Define the functor $D_{loc}=\prod D_{P_i}$. J. Bertin and A. M\'ezard proved that there is a smooth morphism 
$\phi:D_{gl} \rightarrow D_{loc}$, and this morphism induces the following relation on the global deformation ring $R_{gl}$ 
and of the  deformation rings $R_i$ of the deformation functors $D_{P_i}$.
\[
 R_{gl} =(R_1\hat{\otimes} R_2 \hat{\otimes} \cdots \hat{\otimes} R_r)[[U_1,\ldots,U_N]],
\]
where $N=\dim_k H^1(X/G,\pi_*^G(\T_X))$, and $R_i$ is the deformation ring of $D_{P_i}$. For more information concerning this construction we refer to \cite{Be-Me}. For an exact formula for $N$ we refer to \cite[sec. 3]{KontoANT}.

In the approach of Schlessinger \cite{Sch} one wants to build deformations of $(X,G)$
over Artin algebras, especially over the algebras $k[\epsilon]/\epsilon^n$, 
and study whether a deformation over $\Spe k[\epsilon]/\epsilon^n$ can 
be lifted to deformation over $\Spe k[\epsilon]/\epsilon^{n+1}$. More 
generaly 
a small extension $A'$ of $A$ is given by the  
the short exact sequence of local Artin algebras
 \[
0\rightarrow \mathrm{ker}\pi \rightarrow A' \rightarrow A \rightarrow 0
\]
such that  $\mathrm{ker}\pi\cdot m_{A'}=0$, 
where $m_A$ is  the maximal ideal of $A'$ respectively.
We would like to know if a deformation in $D(A)$ can be lifted to 
a deformation in $D(A')$.
The obstructions of such  liftings are elements in $H^2(G,\TO)$.
If there are no obstructions then we can construct a family over 
the formal scheme
$\mathcal{X}\rightarrow \mathrm{Spf} R$ for some complete domain $R$. 
The scheme $\mathrm{Spf} R$ 
is a formal scheme and does not posses a generic fibre. 
J. Bertin and A. M\'ezard in \cite{Be-Me} observed that an 
algebraization theorem of Grothendieck \cite{GroFGA} gives that 
the formal scheme representing $D_{gl}$ is algebraizable, 
and it corresponds to the formal completion of a proper 
smooth curve over $\Spe R$. This means that every 
unobstructed deformation over a formal affine scheme can 
be extended to the generic fibre.
 
Assume that $\mathcal{X} \rightarrow \Spe R$ is a relative 
curve that is a solution to our deformation problem, where $R$ is 
a complete local domain.
Let $\sigma \in G_1(P)$, $\sigma \neq 1$, and let $\tilde{\sigma}$ be a lift 
of $\sigma$ in $\mathcal{X}$. The scheme $\mathcal{X}$ is regular at $P$,
 and the completion of $\mathcal{O}_{\mathcal{X},P}$ is isomorphic to the ring $R[[T]]$.
 Weierstrass preparation theorem \cite[prop. VII.6]{BourbakiComm}   implies that:
\[
\tilde{\sigma}(T)-T=g_{\tilde{\sigma}}(T) u_{\tilde{\sigma}}(T),
\]
where $g_{\tilde{\sigma}}(T)$ is a distinguished Weierstrass polynomial 
of degree $m+1$ and $u_{\tilde{\sigma}}(T)$ is a unit in $R[[T]]$.

The polynomial $g_{\tilde{\sigma}}(T)$ gives rise to a horizontal divisor 
that corresponds to the fixed points of $\tilde{\sigma}$. This 
horizontal divisor might not be reducible. 
The branch divisor corresponds to the union of the fixed points of 
any $\sigma \in G_1(P)$. 
Next lemma shows how to define a horizontal branch divisor for the 
relative curves $\mathcal{X} \rightarrow \mathcal{X}^G$ when 
$G$ is not a cyclic group.

\begin{lemma} \label{lemmaBRANCH}
Let $\mathcal{X} \rightarrow \Spe A$ be an $A$-curve, admitting a 
fibrewise action of the finite group $G$, where $A$ is a 
Noetherian local ring. 
Let $S=\Spe A$, and $\Omega_{\mathcal{X}/S}$, $\Omega_{\mathcal{Y}/S}$ be
the sheaves of relative differentials of $\mathcal{X}$ over $S$ and 
$\mathcal{Y}$ over $S$, respectively. Let $\pi:\mathcal{X} \rightarrow 
\mathcal{Y}$ be the quotient map.
The sheaf 
\[
\mathcal{L}(-D_{\mathcal{X}/\mathcal{Y}})= \Omega_{\mathcal{X}/S} ^{-1}
\otimes_S \pi^* \Omega_{\mathcal{Y}/S}. 
\]
is the ideal sheaf the horizontal Cartier divisor 
 $D_{\mathcal{X}/\mathcal{Y}}$. The intersection of $D_{\mathcal{X}/\mathcal{Y}}$ with the special and generic fibre 
of $\mathcal{X}$ gives the ordinary branch divisors for curves.
\end{lemma}
\begin{proof}
We will first prove that 
 the above defined divisor $D_{\mathcal{X}/\mathcal{Y}}$ is indeed 
an effective Cartier divisor. According to \cite[Cor. 1.1.5.2]{KaMa}
it is enough to prove that 
\begin{itemize}
\item  $D_{\mathcal{X}/\mathcal{Y}}$ is a closed subscheme which is flat over $S$. 
\item for all geometric points $\Spe k \rightarrow S$ of $S$, the 
closed subscheme $D_{\mathcal{X}/\mathcal{Y}}\otimes_S k$ of $\mathcal{X} \otimes_S k$ is a
Cartier divisor in $\mathcal{X} \otimes _S k/k$. 
\end{itemize}

We are interested in deformations of nonsingular curves. 
Since the base is a local ring and the special fibre is nonsingular, 
the deformation $\mathcal{X} \rightarrow \Spe A$ is smooth.  
(See the remark after the definition 3.35 p.142 in \cite{LiuBook}).
The smoothness of the curves $\mathcal{X}\rightarrow S$, 
and $\mathcal{Y}\rightarrow S$, implies that the sheaves 
$\Omega_{\mathcal{X}/S}$ and $\Omega_{\mathcal{X}/S}$ are $S$-flat,
\cite[cor. 2.6 p.222]{LiuBook}. 

On the other hand the sheaf $\Omega_{\mathcal{Y},\Spe A}$ is 
by \cite[Prop. 1.1.5.1]{KaMa}  $\O_{\mathcal{Y}}$-flat. 
Thus, $\pi^*(\Omega_{\mathcal{Y}, \Spe A})$ is $\O_{\mathcal{X}}$-flat
and therefore $\Spe A$-flat \cite[Prop. 9.2]{Hartshorne:77}.
Finally, observe  that  the intersection with the special and generic 
fibre is the ordinary branch divisor for curves according to 
\cite[IV p.301]{Hartshorne:77}.
\end{proof}

{\bf Remark:}
Two horizontal branch divisors can collapse to the same point in 
the special fibre. For instance, this always happens if a 
deformation of  curves from positive characteristic to characteristic zero
with a wild ramification point is possible.

For a curve $X$ and a branch point $P$ of $X$ we will 
denote by  $i_{G,P}$  the order function of the filtration of $G$ at $P$. 
The Artin  representation of the group $G$ is defined 
by $\mathrm{ar}_P(\sigma)=-f_P i_{G,P}(\sigma)$ for $\sigma\neq 1$ and 
$\mathrm{ar}_P(1)= f_P\sum_{\sigma\neq 1} i_{G,P}(\sigma)$ \cite[VI.2]{SeL}.
We are going to use the Artin representation at both the special 
and generic fibre. In the special fibre we always have $f_P=1$ since the
field $k$ is algebraically closed. The field of quotients of $A$ should 
not be algebraically closed therefore a fixed point there might have $f_P \geq 1$.
The integer  $i_{G,P}(\sigma)$ 
is equal to the multiplicity of $P\times P$ in the intersection of 
$\Delta .\Gamma_\sigma$ in the relative $A$-surface 
$\mathcal{X} \times_{\Spe A} \mathcal{X}$, 
where $\Delta$ is the 
diagonal and $\Gamma_\sigma$ is the graph of $\sigma$ \cite[p. 105]{SeL}. 
  
Since the diagonals  $\Delta_0,\Delta_\eta$  and the graphs of $\sigma$ in the special and generic fibres respectively of 
$\mathcal{X}\times_{\Spe A} \mathcal{X}$  are  algebraically equivalent divisors   we have: 
\begin{pro}\label{bertin-gen}
	Assume that $A$ is an integral domain, and let $\mathcal{X}\rightarrow \Spe A$
	be a deformation of $X$. 
        Let $\bar{P}_i$, $i=1,\cdots,s$ be the  horizontal branch divisors 
        that intersect at the special fibre, at point $P$, and let $P_{i}$ be 
        the corresponding points on the generic fibre. For the Artin 
        representations attached to the points $P,P_{i}$ we have:
        \[
        \mathrm{ar}_P(\sigma)=\sum_{i=1}^s \mathrm{ar}_{P_{i}}(\sigma). 
        \]
\end{pro}
This generalizes a result of  J. Bertin \cite{BertinCRAS}. Moreover 
if we set $\sigma=1$ to the above formula we obtain a relation 
for the valuations of the differents 
 in the special and the generic fibre, since the value 
of the Artin's representation at $1$ is the valuation of 
the different \cite[prop. 4.IV,prop. 4.VI]{SeL}.  This observetion 
is equivalent to claim  3.2 in \cite{MatignonGreen98} and is one
 direction of a local 
criterion for good reduction theorem proved in \cite[3.4]{MatignonGreen98},
\cite[sec. 5]{KatoDuke87}.

\begin{cor} \label{artin-lift0}
 Assume that  $V=G_1(P)$ is an elementary abelian group with more than one $\Z/p\Z$ components. 
If $V$ can be  lifted to characteristic zero, 
then $\frac{|V|}{p} \mid m+1$. 
\end{cor}
\begin{proof}
The group $V$ acts on the generic fibre, where the possible stabilizers of points are 
cyclic groups. Since $V$ is not cyclic it can not fix any point $P_i$  in 
 the intersection of the branch locus with the  generic fibre.
Only a cyclic component of $V$ can fix a point $P_i$. Since $V$ act on 
the set of points $P_i$, each orbit has $|V|/p$ elements.
For any element $\sigma \in V$  the Artin representation $\mathrm{ar}_{P_i}(\sigma)=1$
 (no wild ramification at the generic fibre). Therefore proposition 
\ref{bertin-gen} gives us that the number of $\{P_i\}$ is $m+1$ and the desired result follows.
\end{proof}

{\bf Remark:} \label{ex-cor} Consider the case of equicharacteristic deformations of
ordinary curves, together 
with a $p$-subgroup of the group of automorphisms. 
Then $|\mathrm{ar}_P(\sigma)|=2$  for all 
$\sigma \in G(P)=G_1(P), \sigma\neq 1$ \cite{Nak}.
On the other hand the ramification at the 
points of the generic fibre is also wild and \ref{bertin-gen} implies that 
there is only one horizontal branch divisor extending every wild ramification 
point $P$.

{\bf Remark:} The author finds amusing the following similarity 
to the theory of dynamical systems: It is known that autonomous 
(ordinary) differential equations on a manifold $M$ induce an action
of $\mathbb{R}$ on $M$. The fixed locus of this action, called 
{\em equilibrium locus} in the realm of differential equations, can split as the integrated vector fields 
depend on parameters. The study of this splitting  is the object of  {\em
bifurcation theory} \cite{HaleKocak}.
 Notice also that $\mathbb{R}$ is not compact 
and the representation theory of $\mathbb{R}$ shares many difficulties 
with the corresponding representation theory of groups of order divided 
by the characteristic, because of the absence of a Haar measure on them.

\begin{pro} \label{main-free}
Let $R$ be a complete local regular integer domain.
Let $\mathcal{X}\rightarrow \Spe R$ be a deformation of the couple $(X,G)$, 
and let $P$ be a wild ramified point of the special fibre $X$. Assume that 
there is a 
  a $2$-dimensional representation 
$\rho:G_1(P) \rightarrow \mathrm{GL}_k(H^0(X,\mathcal{L}(mP)))$ attached to $P$. Assume also that there is a $G$-invariant horizontal divisor 
that intersects the special fibre with multiplicity $m$.
Then, 
there is a free $R$-module $M$ of rank $2$ generated by $1,\tilde{f}$
so that $M:=\langle 1, \tilde{f} \rangle_R \subset H^0((\mathcal{X},\mathcal{L}(\alpha D))),$
where $1\leq \alpha \in \mathbb{N}$ and $M\otimes_R k=H^0(X,\mathcal{L}(mP))$.
Moreover,
the representation $\rho$ can be lifted to a representation 
\[\tilde{\rho}:G_1(P) \rightarrow \mathrm{GL}_R(   
\langle 1, \tilde{f} \rangle_R). \]
The elements $\tilde{\rho}_\sigma$  are  lower triangular matrices.

Moreover the basis element $\tilde{f}$ is of the form 
\begin{equation} \label{F-form}
\tilde{f}=\frac{1}{(T^m+a_{m-1} T^{m-1}+\cdots + a_1 T_1+a_0)}u(T),
\end{equation}
where $a_0,\ldots,a_{m-1} \in m_R$ and $u(T)$ is a  unit 
in $R[[T]]$ reducing to $1 \mod m_R$.
\end{pro}
\begin{proof}
Let us consider the sheaf $\mathcal{L}(D)$. The space of global sections $H^0(\mathcal{X},\mathcal{L}(D))$
has the structure of  an $R$-module.  For an arbitrary Cartier divisor $D$ on $\mathcal{X}$  and for all $i\geq 0$ there is a natural map \cite[prop. III 12.5]{Hartshorne:77}
\[
\phi_i: H^i(\mathcal{X},\mathcal{L}(D)) \otimes _R k \rightarrow 
H^i(X_s,\mathcal{L}(D\otimes k)).\]  
We are interested in global sections {\em i.e.},  for the  zero  cohomology groups, but in general $\phi_0$ 
can fail to be an isomorphism.

Instead of looking at $D$ we will consider
 $a'D$, where $a'$ is a sufficiently large natural number. 
 We will employ the Riemann-Roch theorem in both the special 
 and the generic fibre and we can choose $a$ sufficiently big 
 so that the index of speciality at both the generic and the special fibre is zero. 
P. Deligne - D. Mumford observed \cite[4. 78]{DelMum},
\cite[chap.3 sec.7]{EGAIII1} that since 
\[H^1(\mathcal{X}_s,\mathcal{L}(a'D \otimes k))=
H^1(\mathcal{X}_\eta,\mathcal{L}(a'D \otimes K))=0\]
the $R$-module $ H^0(\mathcal{X},\mathcal{L}(a'D))$ is free.
We can then select an element $\tilde{f} \in H^0(\mathcal{X},\mathcal{L}(a'D))$
so that $\tilde{f} \equiv f \mod m_R$. 
Consider the least $a$ such that 
$\langle 1,\tilde{f} \rangle_R \subseteq H^0(\mathcal{X}, \mathcal{L}(aD))$ for some 
$1\leq a \leq a'$. 
Since $D$ is $G_1(P)$-invariant the $R$-module 
$H^0(\mathcal{X}, \mathcal{L}(aD))$  is 
equipped with a $G_1(P)$-action.  
The module $M$ might not be the whole 
$H^0(\mathcal{X}, \mathcal{L}(aD))$ but it is the $R$-free part of it.
Therefore $G_1(P)$ acts on $M$ as well and 
 the representation can be lifted:
\[
\tilde{\rho}:G_1(P) \rightarrow \mathrm{GL}_R(M),
\] 
as required.
Since $\sigma\mid_R=\mathrm{Id}_A$ this representation is given 
by lower triangular matrices.

The element $1/\tilde{f}$ is a holomorphic element in $R[[T]]$ reducing to 
$1/f=t^m$ modulo $m_R$. Thus, the reduced order of $1/\tilde{f}$ is $m$ and 
eq. (\ref{F-form}) follows by Weierstrass preparation theorem \cite[prop. VII.6]{BourbakiComm}.
\end{proof}
We will now try to give conditions for the existence of a $G_1(P)$-invariant
divisor intersecting the special fibre at $P$ with degree $m+1$. 
Let $T=\{\bar{P}_i\}_{i=1,\ldots,s}$ be the set of horizontal 
branch divisors that restricts to $P$ in the special fibre of $X$.
This space is acted on by $G_1(P)$, since $\bar{P}_i$ are all 
components of the branch divisor. Each of the $\bar{P}_i$ is fixed 
by some element of $G$ but not necessarily by the whole 
group $G_1(P)$, unless of course $G_1(P)$ is isomorphic to $\Z/p\Z$.

Let $O(T)$ be the set of orbits of $T$ under the action of the group $G_1(P)$,
on $T$. 
A horizontal divisor $D$ supported on $T$, is invariant under the action of 
$G_1(P)$ 
if and only if, the divisor $D$ is of the form:
\begin{equation} \label{pp-oo13}
D=\sum_{C\in O(T)} n_C \sum_{P\in C} P, 
\end{equation}
{ i.e.}, horizontal Cartier divisors  that are in the same orbit of
the action of $G_1(P)$ must appear with the same weight in $D$. 
If the semigroup $\sum_{C\in O(T)} n_C \# C$,  $n_C \in \mathbb{N}$, contains 
        the Weierstrass  semigroup of the branch point $P$ of the special 
fibre, then we can select the desired $G_1(P)$-invariant divisor $D$ supported on 
$T$.
         
If one orbit of $G_1(P)$ acting on $T$ is a singleton, i.e., 
there is a $\bar{P}_i$ fixed by the whole group $G_1(P)$, then
the semigroup
\[\sum_{C\in O(T)} n_C \# C, \;\;\; n_C \in \mathbb{N},\]
is the semigroup of 
natural numbers, and we are done.
This is the case when the group $G_1(P)$ is cyclic.

 If $\#T \not\equiv 0 \mathrm{mod} p$ then there is at least 
one orbit that is a singleton. 
Indeed, 
if all orbits have more than one element then all orbits must have 
cardinality  divisible by $p$, and since the set $T$ is the disjoint union 
of orbits it must also have cardinality divisible by $p$.
\begin{lemma} \label{lem2.5}
If $m$ is the first pole number that is not divisible by the characteristic, and $p\nmid m+1$ then 
there is an orbit that consists of only one element.
\end{lemma}
\begin{proof}
By  proposition \ref{bertin-gen} the 
Artin representation at the special fibre equals the sum of the 
Artin representations at the generic fibre.
Let $\sigma \in G_1(P)$. The Artin representation of $\sigma$ at the 
special fibre equals $m+1$. All $\bar{P}_i$ that are not 
fixed by $\sigma$ do not contribute in the sum of the 
Artin representations at the generic fibre. 

An element $\tau$ sends $\bar{P}_i$ which  is fixed by $H \subset G_1(P)$
to $\tau \bar{P}_i$ which is fixed by $ \tau H \tau^{-1}$. 
Since the representation attached to $P$ is two dimensional the 
group $G_1(P)$ is abelian, and $\tau \bar{P}_i$ is fixed by $H=\tau H \tau^{-1}$. 

If we now consider $P_i$ that is fixed by $\langle \sigma \rangle$ 
then the above argument shows that the orbit of $P_i$ under the 
action of the group $G_1(P)$ has $p^a$ elements $0 \leq a$. If $a=0$ 
then $P_i$ is fixed by the whole group $G_1(P)$. If on the other 
hand for all $P_i$ fixed by $\sigma$ the coresponding 
orbit orders have more than one element then  
the set of $\bar{P}_i$ fixed by $\sigma$ has order divisible by $p$.
This implies that the sum of the 
 Artin representations
at the generic fibre is divisible by $p$, a contradiction.  
\end{proof}
We have thus obtained the following easy to apply 
\begin{cor} \label{cor2.6}
If $G_1(P)$ is cyclic 
or $p\nmid m+1$, then there is a 
horizontal branch divisor $D$, fixed under the action of $G_1(P)$,
 that intersects the special fibre at $mP$. 
In particular, the 
assumption of proposition \ref{main-free} is satisfied and the 
two dimensional representation can be lifted. 
\end{cor}

\begin{lemma} \label{cycONE}
In the mixed characteristic case, 
if the elementary abelian group $G_1(P)$ has more than two cyclic 
components, then there is no  horizontal $G_1(P)$-invariant  divisor $D$
contained in the branch locus and intersecting the special fibre 
at $P$ with multiplicity $m$.
\end{lemma}
\begin{proof}
Since the stabilizers of elements in the 
generic fibre are cyclic groups of order $p$, 
all orbits of elements are divisible by $p$. Therefore, 
a $G_1(P)$-invariant divisor should have degree divisible by $p$. 
This, can not happen since $(m,p)=1$.
\end{proof}

\begin{rem} \label{n1}
{\em
Lemma \ref{cycONE} shows that our method can not 
be used for lifting curves with elementary abelian action 
to characteristic zero. However, M. Matignon 
 proved that such liftings exist \cite{MatManusc}.
}
\end{rem}

We have seen how to relate a deformation of the couple $(X,G)$ to 
a deformation of a matrix representation. Now we will see the 
effect of considering equivalent deformations of couples.

\begin{lemma} \label{equiv2matrix}
Let $\phi$ be a map $\mathcal{O}_{\mathcal{X},P} \rightarrow \mathcal{O}_{\mathcal{X},P}$
making the extensions $\tilde{\rho}_\sigma, \tilde{\rho}'_\sigma$ equivalent.
The corresponding matrix representations are conjugate by a $2 \times 2$ matrix 
of the form $
\begin{pmatrix}
1 & 0 \\
\mu & \lambda
\end{pmatrix}
$
where $\lambda\equiv 1 \mod m_A$ and $\mu \equiv 0 \mod m_A$. 

Conversely, every such matrix gives rise to a map 
$\phi:\mathcal{O}_{\mathcal{X},P} \rightarrow \mathcal{O}_{\mathcal{X},P}$
that reduces to the identity modulo $m_A$.

\end{lemma} 
\begin{proof}
Assume that there is a map $\phi:\mathcal{O}_{\mathcal{X},P} \rightarrow \mathcal{O}_{\mathcal{X},P}$ making the extensions $\tilde{\rho}_\sigma, \tilde{\rho}'_\sigma$ equivalent.
The local-global principle of J.Bertin-A.M\'ezard implies that this map can be extended to a map 
$\phi':\mathcal{X} \rightarrow \mathcal{X}$ that makes the corresponding global 
deformations equivalent.  Let  $\tilde{f}$ be the generator given in proposition \ref{main-free}. Then $\phi'(\tilde{f}) \in H^0(\mathcal{X},\mathcal{L}(aD))$, therefore
$\phi'(\tilde{f}) =\lambda \tilde{f}+\mu$. This means that $\phi$ gives rise to a 
base change in $  H^0(\mathcal{X},\mathcal{L}(a'D))$, 
and two elements in $F(\cdot)$ are equivalent if they are conjugate by a $2\times 2$ matrix of the desired form.

Conversely, assume that we have two equivalent matrix representations 
that are conjugate by a matrix $Q$ of the form $
\begin{pmatrix}
1 & 0 \\
\mu & \lambda
\end{pmatrix}
$
where $\lambda\equiv 1 \mod m_A$ and $\mu \equiv 0 \mod m_A$.  
Then $Q$ sends  
$\tilde{f} \mapsto \lambda \tilde{f} + \mu$, i.e. 
\[
\frac{1}{\phi(T)^m+\sum_{\nu=0}^{m-1} a_\nu \phi(T)^nu }=
\lambda \tilde{f}(T) + \mu.
\]
A solution $\phi(T)$ of this polynomial equation exists by 
using Hensel's lemma. This solution gives rise to the 
desired map $\phi$.
\end{proof}

\section{Deformations of Linear groups}
\label{secmatdef}
\label{MatrixDeformations}

We would like to represent the functor $F$ defined in Eq. (\ref{Fdeformation}).
We will  employ the construction for universal deformation 
rings for matrix representations,  explained by 
B. de Smit and H. W. Lenstra in \cite{SMLE:97}.
Let $H$ be a $p$-group with identity $e$ and let $\rho: H \rightarrow L_n(k)$ be a faithful representation of $H$.
Let $\Lambda[H,n]$ be the commutative $\Lambda$-algebra generated by $X_{ij}^g$ for 
$g\in H, 1\leq j \leq i \leq n$, such that 
\[
X_{ij}^e=\begin{cases} 1 & \mbox{ if } i=j \\ 0  & \mbox{ if } i \neq j \end{cases}
\]
\begin{equation} \label{ostru-123}
X_{ij}^{gh}= \sum_{l=1}^n X_{il}^g X_{lj}^h \mbox{ for } g,h \in H \mbox{ and }
1 \leq i,j \leq n.
\end{equation}
and 
\[
X_{ij}^g=0 \mbox{ for } i<j \mbox{ and for all } g\in H.
\]
We will focus on representations 
on $L_n(A)$. 
For every $\Lambda$-algebra $A$ we have a canonical 
bijection 
\[
\mathrm{Hom}_{\Lambda-\mathrm{Alg}} (\Lambda[H,n],A)  \cong \mathrm{Hom} (H, L_n (A)),
\]
where a $\Lambda$-algebra homomorphism $f:\Lambda[H,n]\rightarrow A$ 
corresponds to the group homomorphism $\rho_f$ that sends 
$g\in H$ to the matrix $(f(X_{ij}^g))$. 
The representation $\rho: H \rightarrow L_n(k)$ corresponds to 
a  homomorphism $\Lambda[H,n]\rightarrow k$. Its kernel is a maximal 
ideal, which we denote by $m_\rho$. We take the completion 
$R(H)$ of $\Lambda[H,n]$ at $m_\rho$. The canonical map 
$\Lambda[H,n]\rightarrow R(H)$, 
gives rise to a map $\rho_{R(H)}:H \rightarrow L_n(R(H))$, 
such that the diagram:
\[
\xymatrix{
H \ar[r]^{\rho_{R(H)\;\;\;\;\;}}  \ar[d]_{=} & L_n(R(H)) \ar[d] \\
H \ar[r]^{\rho} & L_n(k) 
}
\]
is commutative.

We have to  distinguish two cases:

%
%
$\bullet$ The case of equicharacteristic deformations, i.e.,
$R$ is a complete local domain so that $\mathrm{Quot}(R)$ is of characteristic $p$.
Recall that in this case $\Lambda=k$.
Since the generic fibre is of characteristic $p$ we have 
 $X_{22}^g=1$ for all $1\leq i \leq n$.
Moreover, if we fix elements $g_i$ generating 
$H$ as an $\Z/p\Z$-vector space and monomials 
$x_i=X_{21}^{g_i}-c(g_i)$ for each $g_i$ we easily see that $R(H)=k[[x_1,\ldots,x_n]]$.

%
%
$\bullet$ The case of liftings to characteristic zero, i.e. 
$R$ is a complete local domain so that $\mathrm{Quot}(R)$ 
 of characteristic $0$.
Let us again fix elements $x_i,y_i$ for each generator $g_i$ of $H$,
so that $x_i=X_{21}^{g_i}-c(g_i)$, and $y_i=X_{22}^{g_i}-1$.

In this case we have the conditions:
\begin{equation} \label{e1}
\left(X_{22}^g\right)^p=1,
\end{equation} 
\begin{equation} \label{e2} 
X_{21}^g\sum_{\nu=0}^{p-1} \left(X_{22}^g\right)^\nu=0,
\end{equation} 
and the commuting relation:
$(X_{21}^g-X_{21}^h + X_{22}^gX_{21}^h-X_{22}^hX_{21}^g)=0$. 
Observe that $X_{22}^g\neq 1$. Indeed, if $X_{22}^g=1$ then 
eq. (\ref{e2}) will give us that $X_{21}^g=0$ and then the matrix 
is just the identity. Therefore, equations (\ref{e1}) and (\ref{e2})
reduce to the single equation $\sum_{\nu=0}^{p-1} \left(X_{22}^g\right)^\nu=0$.

These conditions 
imply that:
\[
R(H)=\Lambda[[x_1,\ldots,x_n,y_1,\ldots,y_n]]/I,
\]
where $I$ is the ideal 
\[
I:=\left\langle
\sum_{\nu=0}^p (1+y_i)^{\nu-1}, 
y_j(c(g_i)+x_i)-y_i(c(g_j)+x_j)\right\rangle. 
\]

The ring $R(H)$ defined above does not represent the 
deformation functor $F$, since $A$-equivalent deformations may correspond 
to different maps in $\mathrm{Hom}(R(H),A)$.
If $n=2$, i.e., in the case of a two dimensional 
representation, the conjugation action 
given by lemma \ref{equiv2matrix} is easy to handle.

Considering the quotient of  $R(H)$ in positive characteristic, for representations of 
dimension $\geq 3$ is a difficult problem since  the 
 ''trace`` argument
of characteristic zero does not work. (Characters do not distinguish  equivalent representations 
in modular representation theory).

We focus now on the theory of two dimensional representations.
This forces the group $H$ to be elementary abelian.
We compute that  
\begin{equation} \label{expconj}
\begin{pmatrix}
1 & 0 \\
\mu & \lambda
\end{pmatrix}
\begin{pmatrix}
1 & 0 \\
x & y
\end{pmatrix}
\begin{pmatrix}
1 & 0 \\
\mu & \lambda
\end{pmatrix}^{-1}=
\begin{pmatrix}
1 & 0 \\
\mu+\lambda x-y\mu & y
\end{pmatrix}.
\end{equation}
We will consider the effect of the conjugation action given in 
eq. (\ref{expconj}).
The elements $y_i$ remain invariant while the elements 
$x_i \mapsto x_i+ 
\lambda_a c(g_i)+\lambda_ax_i-\mu y_i$, where $\lambda=1+
\lambda_a$,  $\lambda_a,\mu \in m_A$.
If $A=k[\epsilon]/\epsilon^2$, 
then $x_i=x_i+c(g_i)\lambda_a$ since $\lambda_a x_i \in m_A^2=0$.

Let $A$ be an object  in $\mathcal{C}$.
An element in the set $F(A)$ is determined by the 
conjugation equivalence class of a function $f:R(H)\rightarrow A$. 
Such a function should be defined on the generators $x_i,y_j$ of the 
ring $R(H)$.
Since $f(x_i)$ is equivalent to $f(x_i)+f(\lambda_a) c(g_i) \mod m_A^2$, 
if there is a ring representing the functor $F(\cdot)$
then this should be a subring of $R(H)$ and $f(x_i)=0$ for all generators
$x_i$, as one sees by considering $\lambda_a=-x_i/c(g_i)$. 
Therefore the ring representing $F(\cdot)$ is the subring 
of $R(H)$ generated by $y_1,\ldots,y_n$
in the mixed characteristic case and is the zero ring in the 
equicharacteristic case.
This is in accordance to remark \ref{remd}.  

According to remark \ref{n1} the case $n=1$ is 
the only case we can handle using our 
approach in the mixed characteristic situation.

\begin{rem} \label{nosmooth} 
{\em
We consider the subring  $R$ of
$R(H)=R(\mathbb{Z}/p\mathbb{Z})$ generated by $y$.
The ring $R$ is singular.
Indeed, by the infinitesimal lifting property 
\cite[II. exer. 8.6]{Hartshorne:77}, \cite[sec. 1.4]{HarDef}
it is enough to provide a small extension 
$A'\rightarrow A \rightarrow 0$ and a homomorphism 
$h\in \mathrm{Hom}(R,A)$ that does not lift to a 
homomorphism to $\mathrm{Hom}(R,A')$.
Let $m_\Lambda$ be the maximal ideal of $\Lambda$.
Consider the  
natural map
 $\pi:R \rightarrow R/m_\Lambda R=k[[y]]/\langle y^{p-1} \rangle=:A$.
Consider also the ring $A'$ given by 
$k[[y]]/\langle y^p \rangle$. Then 
$A' \rightarrow A$ is a small extension and there is no map $R \rightarrow A'$ 
lifting $\pi$.  
 Indeed, every such  homomorphism  
$R \rightarrow A'$  should factor  through $\mod m_\Lambda$. Therefore 
we obtain a nontrivial homomorphism $A \rightarrow A'$, 
a contradiction.
}
\end{rem}

{\bf Remark:} In \cite{SOS91} T. Sekiguchi, F. Oort, N. Suwa 
introduced the group schemes $\mathcal{G}^{(\lambda)}$
in order to deform the additive group schemes $\mathbb{G}_a$ to 
the multiplicative group schemes $\mathbb{G}_m$ and they were 
able to give a unified Artin-Schreier-Kummer theory \cite{SeSu94},\cite{SeSu95}. 
Many articles devoted to the deformations of automorphism groups 
from positive to zero characteristic are based on this theory, see 
for example  \cite{MatignonGreen98}.

Observe that if $H=\mathbb{Z}/p\mathbb{Z}$, i.e. 
we have an elementary abelian group with just one component, 
then the ring homomorphism 
\[
R(\mathbb{Z}/p\mathbb{Z}) \rightarrow A[[u,1/(\epsilon u+1)]]
\]
sending $y/\epsilon$ to $ u$ gives rise to 
 an injection of $\hat{\mathcal{G}}^{(\lambda)}\rightarrow \Spe R(\mathbb{Z}/p\mathbb{Z})$, where $\lambda=\epsilon$.
Indeed, the diagonal elements $X_{22}^{g}=1+ \epsilon y/\epsilon=1+\epsilon u$
are multiplied as elements in ${\mathcal{G}}^{(\epsilon)}$.

\section{Relation to first order infinitesimal deformations}
\label{SmallExt}
In this section we will relate the  deformation functor of the two dimensional 
representations given in (\ref{Fdeformation}) to the deformation functor 
of actions in formal powerseries rings in (\ref{Bertin-Mezard-functor}). 
The advantage of this approach is that using the two dimensional 
representation we can contract the infinite powerseries representing the 
extended automorphism to a root of a rational function.
Denote by $V$ the elementary abelian group $G_1(P)$.

Assume that a two dimensional representation is 
attached on the wild ramification point $P$.
By using the equation 
\[
\sigma\left( \frac{1}{t^m}\right)=\frac{1}{t^m} +c(\sigma),
\]
we can define the following representation of $V$ to 
automorphisms of formal powerseries rings:
\[
\rho:V \rightarrow \A(k[[t]]),
\]
\[
\sigma \mapsto \rho_\sigma,
\]
where 
\[
\rho_\sigma(t)=\frac{t}{(1+c(\sigma)t^m)^{1/m}}=
t\left( 
1+ \sum_{\nu=1}^\infty \binom{-1/m}{\nu} c(\sigma)^\nu t^{\nu m}
\right).
\] 
Let \[
0\rightarrow \mathrm{ker}\pi \rightarrow A' \rightarrow A \rightarrow 0
\]
be a small extension, i.e. $\mathrm{ker}\pi\cdot m_{A'}=0$, 
where $m_{A'},m_A$ are the maximal ideals of $A,A'$ respectively.
Assume that we have the following data:
A deformation of the two dimensional representation given by 
$C(\sigma)=c(\sigma)+  \delta(\sigma)$, $\lambda(\sigma)=1+
\lambda_1(\sigma)$, where $\delta(\sigma),\lambda_1(\sigma) \in m_{A'}$
and the element  $\tilde{f}$ given in proposition \ref{main-free}  extending  $f$. Write 
$\tilde{f}=f+ \Delta$, for some element $\Delta \in m_{A'}((t))$. 
Then we have:
\[
\tilde{\rho}_\sigma\left( f+ \Delta) \right)=
\lambda(\sigma)(f+\Delta) +c(\sigma)+ \delta(\sigma).
\]
This implies that  ($f=1/t^m$):
\[
\tilde{\rho}_\sigma \left(\frac{1}{t^m} \right)=
\frac{\lambda(\sigma)}{t^m}+c(\sigma)+
\big( \delta(\sigma)+ \lambda(\sigma)\Delta-\tilde{\rho}_\sigma
\Delta \big),
\]
or equivalently:
\begin{equation} \label{induction}
\tilde{\rho}_\sigma(t)= \rho_\sigma(t)+t \left(
\sum_{\nu=0}^\infty \binom{-1/m}{\nu}
\sum_{k=1}^\nu  \binom{\nu}{k} E^k  c(\sigma)^{\nu-k} t^{m\nu}
\right),
\end{equation}
where 
\[
E=\delta(\sigma)+\lambda(\sigma)\Delta -\tilde{\rho}_\sigma \Delta + 
\frac{\lambda_1(\sigma)}{t^m} \in m_{A'}((t)).
\]
Suppose that we can extend $\rho_\sigma(t)$ to a homomorphism 
$\tilde{\rho}_{\sigma,A} \in \A A[[t]]$. 
A  further extension  of $\rho_\sigma$ over $A'$ is 
then given by 
\[
\tilde{\rho}_{\sigma,A'}(t)=\tilde{\rho}_{\sigma,A}(t)+\rho'_\sigma(t),
\]
where $\rho'_\sigma(t)\in \mathrm{ker}\pi[[t]]$.
Since $\Delta \in m_{A'}((t))$ and since $\mathrm{\ker}\pi\cdot m_{A'}=0$
\[
\tilde{\rho}_{\sigma,A'}(\Delta)=\tilde{\rho}_{\sigma,A}(\Delta).
\]
Thus,
equation (\ref{induction}) allows us to compute  
the value of $\tilde{\rho}_{\sigma,A'}(t)$ from the value of 
$\tilde{\rho}_{\sigma,A}(t)$.

\begin{lemma} \label{matreplift}
Let $\tilde{\rho}_{\sigma,A}=\{\tilde{\rho}_{\sigma,A}(t)\}_{\sigma \in V}$ be a representation 
of $V \rightarrow \A A[[t]]$, and consider the corresponding 
element in $F(A)$. If this element in $F(A)$ can be lifted to 
an element in $F(A')$ then $\tilde{\rho}_{\sigma,A}$ can be lifted to a
representation $V \rightarrow \A A'[[t]]$.
\end{lemma}
\begin{proof}
According to \cite[3.2]{Be-Me} every obstruction in lifting a 
representation in $D(A)$  to $D(A')$ is group theoretic. 
Consider extensions of the 
homomorphisms $\tilde{\rho}_{\sigma,A'} \in \A A'[[t]]$
for every $\sigma \in V$. The element
$\tilde{\rho}_{\sigma,A'} \tilde{\rho}_{\tau,A'} 
\tilde{\rho}_{\sigma \tau,A'}^{-1}$ 
is a $2$-cocycle and gives 
rise to a cohomology class in $H^2(V,\TO)$.

In our case observe that if $\lambda_1(\sigma),\delta(\sigma)$ 
are functions $R(V) \rightarrow A'$ and therefore 
satisfy the $2\times 2$ multiplication relations,
then there is no group theoretic obstruction 
in lifting $\tilde{\rho}_{\sigma,A}$ to $\tilde{\rho}_{\sigma,A'}$ 
since a simple computation shows that the lifts defined by 
eq. (\ref{induction}) satisfy the relations  
\[
\tilde{\rho}_{\sigma,A'} \circ \tilde{\rho}_{\tau,A'}=
\tilde{\rho}_{\sigma\tau,A'}.
\]
Therefore any obstruction to lifting $\{\tilde{\rho}_\sigma\}$
reduces to the corresponding obstruction of lifting the 
matrix representation in $F(A)$ to $F(A')$.
\end{proof}

Now we will focus on the small extension $k[\epsilon]/\epsilon^2 
\rightarrow k$, and we will  compute the image of matrix deformations
in $H^1(V,\TO)$. 
The general cocycle in $H^1(V,\TO)$ is given by $d_1(t) \frac{d}{dt}$.
In \cite{KontoANT} the author proved that the map 
\begin{equation}
\label{myiso}
\TO \rightarrow k[[t]]/t^{m+1}\end{equation}
\[ f(t) \frac{d}{dt} \rightarrow f(t)/t^{m+1}\]
is a $V$-equivariant isomorphism. 
%
%
\begin{pro} \label{5.1}
Assume that $P$ is a wild ramified point of $X$ with a two 
dimensional representation attached to it. 
An extension $\tilde{\rho}_\sigma$ gives rise to the following cocyle in 
$H^1(V,\frac{1}{t^{m+1}}k[[t]])$:
\[
 \alpha(\sigma)=\frac{1}{m} 
\left(
\frac{\lambda_1(\sigma)}{t^m} +\lambda_1(\sigma)c(\sigma)
-\delta(\sigma)
+\sum_{\mu=0}^{m-1}\frac{2m-\mu}{m}
\frac{a_{\mu,1}c(\sigma)}{t^{m-\mu}}
\right),
\]
modulo elements in $A[[t]]$.
\end{pro}
\begin{proof}
We will compute the  first order infinitesimal 
deformations of $\rho_\sigma$.
We begin from 
\[
\tilde{\rho}_{\sigma}(f)=\lambda(\sigma)f+
c(\sigma) +  \delta(\sigma) 
+ \lambda(\sigma)\Delta -  \tilde{\rho}_\sigma\Delta.
\]
Set  $E_1=\frac{\delta(\sigma)}{\lambda(\sigma)}+ \Delta-\tilde{\rho}_\sigma
\Delta \frac{1}{\lambda(\sigma)}-c(\sigma) \lambda_1(\sigma)$. 
Then
\[
\tilde{\rho}_{\sigma}\left(\frac{1}{t^m}\right)=\lambda(\sigma)
\frac{1+t^m \frac{c(\sigma)}{\lambda(\sigma)} + t^m  E'}{t^m}.
\]
We compute 
\begin{eqnarray}
\tilde{\rho}_{\sigma}(t) & = &\frac{\lambda(\sigma)^{-\frac{1}{m}}t}
{
\big(1+t^m c(\sigma)\big)^{1/m} 
\big(
1+\frac{ E_1 t^m}{1+c(\sigma)t^m}
\big)^{1/m}
}  \nonumber \\
&= & 
\frac{\lambda(\sigma)^{-\frac{1}{m}}\rho_\sigma(t) }{\big(
1+\frac{ E_1 t^m}{1+c(\sigma)t^m}
\big)^{1/m}}  \nonumber \\
&=& \lambda(\sigma)^{-\frac{1}{m}}(\rho_\sigma(t) -\frac{1}{m} E_1 \rho_\sigma^{m+1}(t) )\mod \epsilon^2 \nonumber \\
&=& \rho_\sigma(t)-\frac{1}{m} E_1 \rho_\sigma^{m+1}(t) 
- \frac{1}{m}\lambda_1(\sigma) \rho_\sigma(t) \mod \epsilon^2 \label{12}.
\end{eqnarray}
We compute that 
\[
\tilde\rho_\sigma\circ \rho_\sigma^{-1}(t)
=\frac{\tilde{\rho}_\sigma(t)}{(1-c(\sigma)\tilde{\rho}_\sigma(t)^m)^{\frac{1}{m}}}.
\]
Since the derivative of the function 
$x\mapsto \frac{x}{(1+Ax^m)^{\frac{1}{m}}}$ is the function 
$x\mapsto (1+Ax^m)^{-\frac{m+1}{m}}$ we compute: 
\begin{eqnarray}
\left.
\frac{d}{d\epsilon} \tilde\rho_\sigma\circ \rho_\sigma^{-1}
\right|_{\epsilon=0} &=& \frac{t^{m+1}}{\rho_\sigma(t)^{m+1}} 
\left. 
\frac{d}{d\epsilon} \tilde{\rho}_\sigma \right|_{\epsilon=0} \\ 
&= & 
-\frac{1}{m} t^{m+1}  \left. E_1 \right|_{\epsilon=0} 
-\frac{1}{m} \lambda_1(\sigma) \frac{t^{m+1}}{\rho_\sigma(t)^m}  \nonumber \\ 
&=& -\frac{1}{m} t^{m+1}  \left. E_1 \right|_{\epsilon=0} 
-\frac{\lambda_1(\sigma)}{m}  \left( t + t^{m+1} c(\sigma) \right).  \label{der1com}
\end{eqnarray}
We will now compute $(1-\lambda(\sigma)^{-1}\tilde{\rho}_\sigma)\Delta$.
 Write $T=t+\epsilon g_{1}(t)\mod\epsilon^{2}A[[t]].$
Write $\tilde{f}=(T^m+\sum_{\mu=0}^{m-1}a_\mu T^\mu)^{-1}u$, 
where $a_\mu=\sum_{\nu\geq 1} a_{\mu,\nu} \epsilon^\nu$.

We compute:
\begin{eqnarray*}
\Delta=\tilde{f}-\frac{1}{t^m} &= &\frac{1}{T^{m}(1+\sum_{\mu=0}^{m-1}a_\mu T^{\mu-m})}-\frac{1}{t^{m}} \\
&=&\frac{1}{T^{m}}\left(
1-\epsilon \sum_{\mu=0}^{m-1}a_{\mu,1} T^{\mu-m}
\right)-\frac{1}{t^{m}}  \mod\epsilon^{2}A[[T]]\\
&=& 
\frac{1-m\epsilon g_1(t)^{m-1}}{t^m}
\left(
1-\epsilon \sum_{\mu=0}^{m-1}a_{\mu,1} T^{\mu-m}
\right)-\frac{1}{t^{m}}  \mod\epsilon^{2}A[[T]]\\
&=& \epsilon mg_{1}(t)^{m-1}/t^{m}-\epsilon \sum_{\mu=0}^{m-1}a_{\mu,1}t^{\mu-2m}  \mod\epsilon^{2}A[[T]].
\end{eqnarray*}
Consider the automorphism $\sigma$ given by $\sigma(t)=t(1+c(\sigma)t^{m})^{-1/m}.$
Observe that \[
\sigma\left(\frac{1}{t^{k}}\right)=\frac{(1+c(\sigma)t^{m})^{\frac{k}{m}}}{t^{k}}=\frac{1}{t^{k}}+\sum_{\nu\geq1}\binom{\frac{k}{m}}{\nu}c(\sigma)^{\nu}t^{m\nu-k},\]
therefore\[
(1-\lambda_1(\sigma)\epsilon)\sigma\left(\frac{1}{t^{k}}
\right)-\frac{1}{t^{k}}=\frac{k}{m}c(\sigma)t^{m-k}+\sum_{\nu\geq2}\binom{\frac{k}{m}}{\nu}c(\sigma)^{\nu}t^{m\nu-k}-\frac{\epsilon \lambda_1(\sigma)}{t^k}.\]
This means that for $k\leq m$ the quantity $(1-\lambda(\sigma)^{-1}\tilde{\rho}_\sigma)(\epsilon t^{-k})$
is holomorphic in $t$ modulo $\epsilon^2$. Thus $(1-\lambda(\sigma)^{-1}\tilde{\rho}_\sigma)\frac{g_{1}(t)^{m-1}}{t^{m}}\in\mod A[[t]]$
and we arrive at:
\[
(\lambda(\sigma)^{-1}\sigma-1)\epsilon \Delta=\sum_{\mu=0}^{m-1}\frac{2m-\mu}{m}
\frac{a_{\mu,1}c(\sigma)}{t^{m-\mu}}
\mod \epsilon^2 +A[[t]].
\]
This result combined with eq. (\ref{der1com}) gives us
\begin{equation} \label{derfr2com}
\left.
\frac{d}{d\epsilon} \tilde{\rho}_\sigma\circ \rho_\sigma^{-1}
\right|_{\epsilon=0}=
\frac{t^{m+1}}{m} 
\left(
\frac{\lambda_1(\sigma)}{t^m} +\lambda_1(\sigma)c(\sigma)
-\delta(\sigma)
+\sum_{\mu=0}^{m-1}\frac{2m-\mu}{m}
\frac{a_{\mu,1}c(\sigma)}{t^{m-\mu}}
\right)
\end{equation}
modulo elements in $A[[t]]$.
The desired result follows by applying the map given in eq.(\ref{myiso}).
\end{proof}

%
%
\begin{lemma} \label{4.3}
Assume that $G_1(P)=\mathbb{Z}/p\mathbb{Z}$. The $k$-vector space  $H^1(\mathbb{Z}/p\mathbb{Z},k[[t]]/t^{m+1})$  is generated  by
the elements $\{b_i/t^i: b\leq i \leq m+1\}$ so that $\binom{i/m}{p-1}=0$ and $b=1$ if $p\mid m+1$ and $b=2$ if $p\nmid m+1$, and 
$b_i\in \mathrm{Hom}(\mathbb{Z}/p\mathbb{Z},k)$.
\end{lemma}
\begin{proof}
This is proposition 2.7 in \cite{KontoANT} for $a=-m-1$.
\end{proof}

Consider the elementary abelian group $V=\oplus_{i=1}^s V_i$ where 
$V_i\cong \mathbb{Z}/p\mathbb{Z}$.
The computation of the cohomology group $H^1(V,\TO)$ seems complicated in 
the general case. However, under some mild assumptions we can prove the following:
\begin{pro} \label{cohomologySplit}
Let $m+1=\sum_{i\geq 0} b_i p^i$ be the $p$-adic expansion of $m$.
If $\lf \frac{2 b_0}{p} \rf=\lf \frac{b_0+b_{\nu-1}}{p} \rf$ for all 
$2\leq \nu \leq s$, then the map
\begin{equation} \label{ontomap}
\Psi:H^1(V,\TO) \rightarrow \bigoplus_{\nu=1}^s H^1(V_\nu,\TO),
\end{equation}
sending $v \mapsto \sum_{\nu=1}^s \mathrm{res}_{V \rightarrow V_i}v$
is an isomorphism. Moreover 
\begin{equation} \label{eeqq}
H^1(V,\TO)\cong \bigoplus_{i=2, \binom{i/m}{p-1}=0}^{m+1} b_i \frac{1}{t^i},
\end{equation}
 where   $b_i \in \mathrm{Hom}(V,k)$. 
\end{pro}
\begin{proof}
Consider the maps $c_i\in \mathrm{Hom}(V_i,k)$ and 
extend them to maps $\bar{c}_i\in \mathrm{Hom}(V,k)$, by setting 
$\bar{c_i}(\sigma)=0$ if $\sigma \not\in V_i$. The image of
$\sum_{\nu=1}^s \bar{c}_i$ under the map $\Psi$ 
 given in (\ref{ontomap}) is $(c_1,\ldots,c_s)$,  therefore the map $\Psi$ is  onto
 and it is sufficient to prove that 
both spaces have the same dimension.

For the dimension $h_1(V,\TO)=\dim_k H^1(V,\TO)$ 
the author has proved the following formula:
\begin{equation}\label{eq-el-ab}
h_1(V,\TO)=\sum_{i=1}^s \left( \lf \frac{(m+1)(p-1)+a_i}{p} \rf
-\lc \frac{a_i}{p} \rc \right),  
\end{equation}
where $a_1=-(m+1)$, $a_i=\lc \frac{a_{i-1}}{p} \rc$ \cite[prop. 2.9]{KontoANT}.
Observe that $a_i=-\lf \frac{m+1}{p^{i-1}} \rf$.
We compute that 
\begin{equation} \label{pad1}
\frac{m+1}{p^k}=\sum_{\nu=0}^{k-1} \frac{b_i}{p^{k-\nu}}+\sum_{\nu \geq k}
b_\nu p^{\nu-k}, 
\end{equation}
therefore 
\begin{equation} \label{pad2}
\lf \frac{m+1}{p^k} \rf= \sum_{\nu \geq k}
b_\nu p^{\nu-k}.
\end{equation}
Now we compute that
\begin{equation} \label{pad3}
\lf 
\frac{m+1}{p}+\frac{1}{p} \lf \frac{m+1}{p^{i-1}} \rf
\rf=
\lf\frac{b_0+b_{i-1}}{p}
\rf+
\sum_{\nu \geq 1} b_\nu p^{\nu-1}+\sum_{\nu\geq i} b_\nu p^{\nu-i}.
\end{equation}
The desired result follows by plugging eq. (\ref{pad2}),(\ref{pad3})
into eq. (\ref{eq-el-ab}).

Equation (\ref{eeqq}) folows by lemma \ref{4.3}.
\end{proof}

{\bf Remark:}
Consider the curves defined by
\[
\sum_{\nu=0}^s a_n y^{p^n}=\sum_{\mu=0}^mb_\mu x^\mu,
\]
so that $m\not\equiv 0 \mod p$, $a_s,a_0,b_0\neq 0$, 
$s\geq 1$, $mu\geq 2$ studied by H. Stichtenoth in \cite{StiII}. 
The representation attached to the unique place
$P_\infty$ above the place $p_\infty$ of the function 
field $k(x)$ is two dimensional if and only if  $m< p^s$ \cite{KontoMathZ}.
In this case the assumptions of proposition \ref{cohomologySplit}
hold. 

\begin{cor} \label{obdt}
If $G_1(P)=\mathbb{Z}/p\mathbb{Z}$ or if the assumptions of proposition \ref{cohomologySplit} hold then the tangent vector 
corresponding to $0\neq \frac{d}{dt} \in H^1(V,\TO)$ is an obstructed deformation. \end{cor}
\begin{proof}
The element $\frac{d}{dt}$ corresponds to $\frac{1}{t^{m+1}}\in H^1(V,
\frac{1}{t^{m+1}}k[[t]])$. 
Using proposition \ref{cohomologySplit} we see that it is 
impossible to obtain a vector in the direction of $\frac{1}{t^{m+1}}$
using a matrix representation, i.e. an element in $F(\cdot)$.

 Notice that since we have assumed 
that the representation attached to $P$ is two dimensional 
we have that $m>1$.
\end{proof}

\begin{cor} \label{obell-ab}
Assume that $p\nmid m+1$ and the assumptions of 
proposition \ref{cohomologySplit} hold. 
Assume also that $V$ is an elementary abelian group 
with more than one component.
Using the notation of eq. (\ref{eeqq})
unubstructed deformations should satisfy $b_i(\sigma)=\lambda_i c(\sigma)$
for some element $\lambda_i \in k$.
\end{cor} 
\begin{proof}
Condition $p\nmid m+1$ implies that every deformation is coming from a 
matrix representation \ref{cor2.6} and condition follows by using proposition \ref{5.1}.
\end{proof}

\begin{rem} \label{remd}
{\em 
We see that the 
data $\delta(\sigma)$   of the matrix representation deformation 
do not affect the corresponding element in $H^1(V,\TO)$ since 
they appear as coefficients of $t^0$ in  the cocylce expression of
proposition \ref{5.1} and are cohomologous to zero.
What seems to affect the tangent elements is the coefficients 
of the distinguished Weierstrass  polynomial of the function $\tilde{f}$ defined in 
eq. (\ref{F-form}).

On the other hand
in the case of liftings from characteristic $p$ to characteristic zero
the diagonal element $\lambda_1$ appears as coefficient of the element 
$t\frac{d}{dt}$. This construction is similar to the one
of J.Bertin and A. M\'ezard  \cite[lemme 4.2.2]{Be-Me}.
}
\end{rem}

Following  \cite[th. 4.2.8]{Be-Me} we can prove:
\begin{pro} \label{generalizeBeMe}
If $R_P$ denotes the versal deformation ring at $P$, then there is a surjection 
\begin{equation} \label{Rprime}
R_P \rightarrow W(k)[[y]] 
\left/\left\langle \sum_{\nu=1}^p \binom{p}{\nu}y^{\nu-1}
\right\rangle
\right. :=R'
\end{equation}
The ring $R_P$ is not smooth.
\end{pro}
\begin{proof}
We are in the mixed characteristic case so $V= \langle \sigma \rangle$.
According to 
section \ref{secmatdef}, the ring $R'$
gives rise to  deformation of the two dimensional representation given by 
\[
\tilde{\rho}_\sigma=
\begin{pmatrix}
1 & 0 \\
0 &  1+y
\end{pmatrix}
\]
which in turn gives rise to the deformation
\[
\tilde{\rho}_{\sigma}(t)=\frac{(1+y)^{-\frac{1}{m}}}
{
\left(1+\frac{Et^m}{1+c(\sigma)t^m}\right)^{1/m}} \rho_\sigma(t), 
\]
for a suitable element $E$.
The map $\mathrm{Hom}(R_P,\cdot)\rightarrow D(\cdot)$ is smooth
(in the sence of Schlessinger \cite[def. 2.2]{Sch},
 \cite[p. 278]{MazDef}), therefore 
 there is a map $\phi:R_P \rightarrow R'$.
In order to prove that $R_P$ is not a smooth ring we proceed as follows:
Consider the natural map $\pi:W(k) \rightarrow W(k)/p=k$.
We obtain the following map 
\[
\phi\circ \pi :R_P \rightarrow k[[y]]/\langle y^{p-1} \rangle:=A.
\]
 Consider the ring $A'=k[[y]]/\langle y^p \rangle$.
Then $A' \rightarrow A$ is a small 
extension, and there  is no map $R_P \rightarrow A'$  extending $R_P \rightarrow R' \stackrel{\mod p }{\longrightarrow} A$ by remark \ref{nosmooth}.
In this way we obtain an obstruction to the infinitesimal 
affine lifting for the 
affine scheme $\Spe R_P$ therefore $R_P$ is not smooth.

Alternatively one can compute the obstruction 
as an element in $H^2(V,\TO)$ following \cite[lemme 4.2.3]{Be-Me}.
\end{proof}
\begin{pro}
Assume that the hypotheses  of proposition  \ref{cohomologySplit} hold.
Consider the ring $R_1$ defined by 
\[
R_1=\left\{ 
\begin{array}{ll}
k & \mbox{ in the equicharacteristic case } \\
R' & \mbox{ in the mixed characteristic case (see eq. (\ref{Rprime}))}  
\end{array}
\right.
\]
Let $b=1$ if $p\mid m+1$ and $b=2$ if $p\nmid m+1$.
Let $\Sigma$ be the subset of numbers 
$b\leq i \leq m$ so that $\binom{\frac{i}{m}}{p-1}=0$.
Consider the ring $\bar{R}:= R_1[[X_i:i\in \Sigma]]$ and the 
$k$-vector space $W \subset H^1(V,\TO)/\langle d/dt \rangle$ generated 
by elements $\lambda_i c(\sigma) t^{m+1-i}\frac{d}{dt}$.

There is a surjection $R_P\rightarrow \bar{R}$ 
that induces an isomorphism $W\cong \mathrm{Hom}(\bar{R},k[\epsilon]/\epsilon^2)$.
The Krull dimension of $R_P$ is equal to $\#\Sigma$.
\end{pro}

\begin{proof}
We have observed in corollary  \ref{obdt}  that deformations in the 
direction of $d/dt$ are not coming from matrix representations.
The elements  $\frac{1}{t^i}$ for $i\in \Sigma$ are elements in   $H^1\left(V, \frac{1}{t^{m+1}} k[[t]]\right)$ 
that give rise to elements  $ t^{m+1-i} \frac{d}{dt} \in H^1(V,\TO)$. Every deformation on these directions is 
unobstructed by lemma \ref{matreplift}.
\end{proof}

\section{Relation to Deformations of Artin-Schreier curves}
\label{PRIES}
Let $P$ be a wild ramified point of the cover $\pi:X\rightarrow Y=X/G$
so that the corresponding representation is two dimensional.
In this section we will examine the dependence of the 
Artin-Schreier extension $X \rightarrow X/G_1(P)$ on the 
form of matrix representation $\rho:G_1(P) \rightarrow GL_2(k)$.
Then we will restrict to the germs 
$\mathcal{O}_{X,P} \rightarrow \mathcal{O}_{Y,\pi(P)}$, and we will 
study the relation to the deformation functor introduced in \cite{Pries:04}
by R. Pries.  
The approach of R. Pries is to work with germs of 
curves and to deform the defining Artin-Schreier equation.
Since the germs are  living in local rings, that have only one 
maximal ideal, the effect of splitting the branch locus can not 
be studied. Therefore R. Pries 
considers only deformations that do not split the branch locus. 
According to proposition \ref{bertin-gen} it is impossible
 to lift a wild ramified action to characteristic zero, 
without splitting the branch locus. 
We will now  restrict ourselves 
to the equicharacteristic deformation case.

Let $X$ be  a curve that has a $2$-dimensional representation attached at 
a wild ramified point $P$. Denote by $\{1,f\}$ a basis of the 2-dimensional vector space 
$L(mP)$ where  $m:=v_P(f)$ is the highest jump in the upper ramification filtration.

We would like to write down an algebraic equation for the cover $X \rightarrow X/G_1(P)$.
The representation $c=c_1:G_1(P)\rightarrow k$ is a faithful homomorphism of additive groups.
We consider the action of $G_1(P)$ on $f$:
Let $\Phi(Y)$ be the  additive polynomial with set of  roots $\{c_1(\sigma):\sigma  \in G_1(P)\}$. The polynomial $\Phi(Y)$ can be computed as follows:
The group  $G_1(P)$ is by \cite[sec. 3]{KontoMathZ} elementary abelian so 
we express 
 $G_1(P)$ as an $\mathbb{F}_p$ vector space with basis $\{\sigma_i\}$ such  that 
$G_1(P)=\bigoplus_{i=1}^s \sigma_i \mathbb{F}_p$.

 Let $\Delta(x_{1},\ldots,x_{n})$
denote the Moore determinant:\[
\Delta(x_{1},\ldots,x_{n})=\det\left(\begin{array}{cccc}
x_{1} & x_{2} & \cdots & x_{n}\\
x_{1}^{p} & x_{2}^{p} & \cdots & x_{n}^{p}\\
\vdots &  &  & \vdots\\
x_{1}^{p^{n-1}} & x_{2}^{p^{n-1}} & \cdots & x_{n}^{p^{n-1}}\end{array}\right).\]
The additive polynomial $\Phi$ can be expressed in terms of the Moore determinant:
\[
 \Phi(Y)=\frac{\Delta(c(\sigma_1),\ldots,c(\sigma_s),Y)}{\Delta(c(\sigma_1),\ldots,c(\sigma_s))},
\]
see \cite[lemma 1.3.6]{GossBook},\cite[eq. 3.6]{Elkies99}.
Thus, the cover $X \rightarrow X/G_1(P)$ is given in terms of the generalized Artin-Schreier equation 
\[
 \Phi(Y)=\prod_{\sigma \in G_1(P)} \sigma f=N_{G_1(P)}(f).
\]
We would like to represent the curve as a fibre product of Artin-Schreier curves and then using  Garcia's-Stichtenoth's 
normalization \cite{GarciaSticht1991} to write the curve in the form $y^{p^s}-y=u$, where $u$ is an element in the function field of the 
curve $X/G_1(P)$.

There are elements $y_j\in k(X)$ so that $\sigma_i(y_j)=y_j + \delta_{ij}$. Using this notation we can see that the function field 
$k(X)$ can be recovered as the function field of the  fibre product of the curves $y_i^p-y_i=u_i$.
The constant elements $u_i$ can be computed from the map $c:G_1(P) \rightarrow k$ as follows:
Let $V_i=\bigoplus_{\nu=1,\nu\neq i}^s \sigma_\nu \mathbb{F}_p$.
We compute an  additive polynomial  $\mathrm{ad}_i(Y)$ with roots the $\mathbb{F}_p$-vector space $V_i$ using the Moore determinant:
\[
 \mathrm{ad}_i(Y):=\frac{\Delta(c(\sigma_1),\ldots, \widehat{c(\sigma_i)},\ldots,c(\sigma_s),Y)}{\Delta(c(\sigma_1),\ldots, \widehat{c(\sigma_i)},\ldots,c(\sigma_s))}.
\]
These polynomials are invariants of the curve  and the map $c$. Moreover we compute that $y_i':=\prod_{\sigma\in V_i} \sigma f
=\prod_{v\in V_i} (f-v)=\mathrm{ad}_i(f)$. The element $y_i'$ is invariant under the action of $V_i$ and 
$\sigma_i(y_i')=y_i' +\mathrm{ad}(c(\sigma_i))$. We can normalize by setting $y_i=y_i'/\mathrm{ad}(c(\sigma_i))$. Then, 
\[
 \sigma_j(y_i)=y_i+\delta_{ij}.
\]
Following \cite{GarciaSticht1991}  we choose an $\F_p$ basis $\mu_1,\ldots,\mu_s$  of $\F_{p^s}$ and we set 
$y=\sum_{i=1}^s \mu_i y_i$. We observe that 
the function field can be recovered as the following extension of the field $k(X)^{G_1(P)}$:
\[
 y^{p^s} -y =   N_{G_1(P)} \left(\sum_{i=1}^s  \frac{\mu_i \mathrm{ad}_i(f)}{\mathrm{ad}_i(c(\sigma_i))}  \right)=:u.
\]
The element $u\in k(X)^{G_1(P)}$ is an invariant of the action of $G_1(P)$ on $k(X)$.
Observe that 
\begin{equation} \label{trofod-det}
\Delta\big(c(\sigma_1),\ldots,\widehat{c(\sigma_i)},c(\sigma_s),c(\sigma_i)\big)=(-1)^{s-i}\Delta(c(\sigma_1),\ldots,c(\sigma_s)).
\end{equation}
Let $D$ be the operator sending $x \mapsto x^{p^s}-x$. Since $\mu_i \in \mathbb{F}_{p^s}$ we have  $D(\mu_i x)=\mu_i D(x)$. 
The element $u$ can thus also  be expressed by 
\begin{equation} \label{u1}
 u=\sum_{i=1}^s  \mu_i D \left(\frac{ \mathrm{ad}_i(f)}{\mathrm{ad}_i(c(\sigma_i))} \right)=
\sum_{i=1}^s \mu_i (-1)^{s-i}  D\left(\frac{ \Delta \big(  c (\sigma_1),\ldots,\widehat{c (\sigma_i)},\ldots,c(\sigma_s),f  \big)}
{
\Delta \big(  c (\sigma_1),\ldots,c(\sigma_s) \big)
}
 \right)
\end{equation}
Equation (\ref{u1}) allows us to  express $u$ in terms of the following 
determinant:
\begin{equation} \label{u1det}
u_1=\frac{1}{\Delta \big(  c (\sigma_1),\ldots,c(\sigma_s) \big)}
\det
\begin{pmatrix}
\mu_1 & \mu_2 & \cdots & \mu_s & 0 \\
c(\sigma_1) & c(\sigma_2)& \cdots & c(\sigma_s) & f \\
c(\sigma_1)^p & c(\sigma_2)^p& \cdots & c(\sigma_s)^p & f^p \\
\vdots  & \vdots & & \vdots \\
c(\sigma_1)^{p^{s-1}} & c(\sigma_2)^{p^{s-1}}& \cdots & c(\sigma_s)^{p^{s-1}} 
& f^{p^{s-1}}  
\end{pmatrix},
\end{equation}
\[
u=D(u_1).
\]
%
%
Notice that $u_1$ is a polynomial of $f$ of the form 
\[
u_1(f)=\sum_{\nu=1}^{s} o_\nu f^{p^{\nu-1}},
\]
where $o_\nu$ can be computed, in terms of the function $c$, as minor
determinants of the above matrix. 
Then $u(f)$ is a polynomial of $f$ of the form 
\[
u(f)=\sum_{\nu=1}^{2s} a_i f^{p^{\nu-1}}, 
\]
where for $1\leq \nu$, $a_{\nu+s}=-a_\nu^{p^s}$.

Now consider the relative situation:
Consider the element  $\tilde{f} \in A[[t]][t^{-1}]$ defined in 
proposition \ref{main-free}. 
Given such an element $\tilde{f}$ and a deformation of the representation $\rho:G_1(P) \rightarrow GL_2( L(mP))$ 
we will construct a deformation   $\mathcal{O}_{X,P}$  of the germ $\mathcal{O}_{X,P}$ with Galois group $G_1(P)$. 

We  form again the additive polynomials:
\[
 \mathrm{Ad}_i(Y)=: \frac{\Delta(C(\sigma_1),\ldots, \widehat{C(\sigma_i)},\ldots, C(\sigma_s),Y) }{\Delta(C(\sigma_1),\ldots, \widehat{C(\sigma_i)},\ldots, C(\sigma_s))}.
\]
Using the previous normalization procedure  we arrive at the following deformed Artin-Schreier curve:
\begin{eqnarray*}
y^{p^s} -y  &=&   \sum_{i=1}^s \mu_i D\left( \frac{ \mathrm{Ad}_i(\tilde{f})}{\mathrm{Ad}_i(C(\sigma_i))} \right) =\\
                  &=& \sum_{i=1}^s \mu_i (-1)^{s-i} D \left( \frac{ \Delta \big(  C(\sigma_1),\ldots,\widehat{C (\sigma_i)},\ldots,C(\sigma_s),\tilde{f} \big)} 
{
\Delta \big(  C (\sigma_1),\ldots,C(\sigma_s) \big)
}\right)
 := U.
\end{eqnarray*}
Notice that similar to equation (\ref{u1det}) we have:
\[
U_1=\frac{1}{\Delta \big(  C (\sigma_1),\ldots,C(\sigma_s) \big)}
\det
\begin{pmatrix}
\mu_1 & \mu_2 & \cdots & \mu_s & 0 \\
C(\sigma_1) & C(\sigma_2)& \cdots & C(\sigma_s) & \tilde{f} \\
C(\sigma_1)^p & C(\sigma_2)^p& \cdots & C(\sigma_s)^p & \tilde{f}^p \\
\vdots  & \vdots & & \vdots \\
C(\sigma_1)^{p^{s-1}} & C(\sigma_2)^{p^{s-1}}& \cdots & C(\sigma_s)^{p^{s-1}} 
& \tilde{f}^{p^{s-1}}  
\end{pmatrix},
\]
\[
U=D(U_1).
\]
The element $U \in A[[t]][t^{-1}]$ so that $U \equiv u \mod m_A$.

\subsection{Relation with equivalence class of Artin-Schreier extensions.}
In what follows we would like to consider isomorphism  classes of Artin-Schreier curves. The following 
lemma identifies two Artin-Schreier extensions of the ring $A[[x]][x^{-1}]$, where $A$ is a $k$-algebra 
that gives rise to an irreducible affine scheme, i.e. $A/\mathrm{rad}(A)$ is an integral domain.
\begin{lemma} \label{RPlemma}
Consider the extensions $y_1^{p^s}-y_1 =g_1$ and $y_2^{p^s}-y_2=g_2$, 
where $g_1,g_2 \in A[[x]][x^{-1}]$.
 These extensions are isomorphic  if and only if 
$g_1(x)=\zeta g_2(x)+ d^{p^s}-d$, for  some $d \in A[[x]][x^{-1}]$, and $\zeta \in \mathbb{F}_{p^s}^*$.
\end{lemma}
\begin{proof}
 If $A$ is a field  $k$ then this is classical result due to Hasse \cite{Hasse34}. For the general case we refer to
\cite[lemma 2.4]{Pries:05}.
\end{proof}
 \begin{lemma}
 The Artin-Schreier curve $y^{p^s}-y=f(x)$ where $f(x) \in A[[x]][x^{-1}]$ is isomorphic to $y^{p^s}-y=f(x)+g(x)$, where $g(x) \in A[[x]]$.
\end{lemma}
\begin{proof}
 Following \cite[sec. 3]{Pries:05} we observe that $g(x)=d^{p^s}-d$, where $d=\sum_{\nu=0}^\infty g(x)^{p^{s\nu}}$. The desired result 
follows by using lemma \ref{RPlemma}.
\end{proof}
%
%
Let $m$ be the conductor, i.e. the highest jump in the upper ramification filtration. Since the group $H$ is 
elementary abelian this is equal to the highest jump in the lower 
ramification
filtration \cite[lemma 1.8]{KontoANT}.
\begin{lemma} \label{holnotalter}
 Consider an Artin-Schreier cover of $A[[x]][x^{-1}]$ given by:
\[
 y^{p^s}-y =\sum_{\nu=0}^\lambda r_\nu(1/x)^{p^\nu},
\]
where $r_\nu(T)\in A[T]$ are polynomials of degree $d_\nu$, so that $\gcd(d_\nu,p)=1$ . The conductor of the cover equals to 
$\max_{\nu} {d_\nu}$.
\end{lemma}
\begin{proof}
 R. Pries \cite{Pries:05}.
\end{proof}

%
%
D. Harbater in   \cite{Harbater:80} (see also \cite[sec. 5.1]{Be-Me}) gave a parametrization of the  classes of cyclic $\mathbb{Z}_p$-covers
of a local fields branched above the maximal ideal. 
For the more general case of $\mathbb{F}_{p^s}$-covers 
the space of classes of covers of $k((t'))$  is parametrized by the quotient:
\begin{equation} \label{classesHarbater}
 C=\frac{k((t'))}{k[[t]]+D(k((t')))},
\end{equation}
where $D$ denotes the map $x \mapsto x^{p^s}-x$.
Indeed, by lemma \ref{RPlemma} adding $D(a)$ does not alter the equivalence class of the Artin-Schreier curve and 
by lemma \ref{holnotalter} the same holds for adding a holomorphic element.

%
%
R. Pries  gave  a moduli interpretation of $p$-group covers of the projective line and 
she proposed two approaches: either 
 transform (by extending the base ring $A$) an arbitrary Artin-Schreier 
extension of $A[[x]][x^{-1}]$ to a class in  (\ref{classesHarbater}) or 
 define 
a fine moduli space by considering a category where all powers of the $q$-Frobenious maps are 
invertible elements.
She  introduced the following:
\begin{orism}
 Let $A_1,A_2$ be two $k$-algebras that give rise to irreducible affine schemes, i.e. $A_i/\mathrm{rad}(A_i),$ $ i=1,2$ are  
integral domains. Consider the Artin-Schreier relative $A$-curves
$C_i:y_i^{p^s}-y_i=f_i(x)$, where $f_i(x) \in A_i$. The two curves are considered to be equivalent if and only if there is 
an  algebra extension $A$ of both $A_i$, i.e. there are ring monomorphisms  $A_i \hookrightarrow A$, so that 
the curves $C_i \times_{\Spe A_i} \Spe A$ are isomorphic covers of $A[[x]][x^{-1}]$. 
\end{orism}

In 
general $U-u \in m_A [[t]][t^{-1}]$ and it is not 
an element in $m_A [[x]][x^{-1}]$. 
 If the deformation $\tilde{\rho}_\sigma$
does not split the branch locus, then $U-u \in m_A [[x]][x^{-1}]$.
 After cutting the holomorphic 
part of $U-u$ and applying the transformation of lemma \ref{RPlemma} we get an 
equivalence class of germs of Artin-Schreier curves given in eq. 
(\ref{classesHarbater})  to an element  in the deformation functor of 
Pries.



Conversely for every Laurent polynomial $\Delta \in  m_A((x))$
so that $n_0=v_x(\Delta)$, satisfies $(n_0,p)<m$ we can consider the 
extension of $A((x))$ defined as 
\[
A((x))[y]/(y^{p^s}-y=f+\Delta).
\] 
This gives rise to an infinitesimal extension of the germ of $X$ 
at $P$ in the sense of Pries and according to the local-global 
theory developed  by Harbater all this local deformations 
can be patched together to give a global deformation of the 
couple $(X,G)$.

\def\cprime{$'$}
\providecommand{\bysame}{\leavevmode\hbox to3em{\hrulefill}\thinspace}
\providecommand{\MR}{\relax\ifhmode\unskip\space\fi MR }
\providecommand{\MRhref}[2]{%
  \href{http://www.ams.org/mathscinet-getitem?mr=#1}{#2}
}
\providecommand{\href}[2]{#2}

\end{document}